[1994/12/01]
\documentclass[draft]{article}
\title{Stochastic impulsive fractional differential evolution equations with infinite delay }
\author{Shufen Zhao   Minghui Song}

\usepackage{amsmath,amsthm}

\chardef\bslash=`\\ 

\newtheorem{thm}{Theorem}[section]

\newtheorem{lem}[thm]{Lemma}

\theoremstyle{definition}
\newtheorem{defn}{Definition}[section]

\theoremstyle{remark}
\newtheorem{rem}{Remark}[section]

\newcommand{\eval}[2][\right]{\relax
  \ifx#1\right\relax \left.\fi#2#1\rvert}

\begin{document}
\maketitle
\markboth{Stochastic impulsive fractional differential evolution equations ...}
{Stochastic impulsive fractional differential evolution equations ...}
\renewcommand{\sectionmark}[1]{}
\begin{abstract}
In this paper, we investigate a class of stochastic impulsive fractional differential evolution equations with infinite delay in Banach space. Firstly sufficient conditions of the existence and uniqueness of the mild solution for this type of equations are derived by means of the successive approximation.  Then we use the Bihari's inequality to get the stability in mean square of the mild solution. Finally an example is presented to illustrate the results.

{\bf Mathematics Subject Classfication} 34K30; 34F05; 60H30.
\end{abstract}
\section{Introduction}

In resent years, the differential equations of fractional order have been widely studied by many authors (see \cite{a,b,3,c,d,e} and references therein ) due to their application in many practical dynamical phenomena arising in engineering, physics, economy and science \cite{f,g,h,i,J}. Very recently, the author of \cite{12} established the sufficient conditions for the existence and uniqueness of mild solution for impulsive fractional integro-differential evolution equations with infinite delay by means of the Kuratowski measure of noncompactness and progressive estimation method.

It is well known that time delay phenomena are frequently encountered in a variety of dynamic systems such as nuclear reactors, chemical
engineering systems, biological systems and population dynamic models, at the same time environment noise often perturbed many branches of science, for
this reason when
we take environment noise and time delays into account (see \cite{liuk,tt,8,hale, hou,16,cc,bb, 7}), it is reasonable to consider the stochastic fractional evolution equations with delays.

This paper is concerned with the existence and uniqueness of mild solutions for Cauchy problems for the stochastic impulsive fractional evolution equations:
\begin{eqnarray}\label{eq1}
\left\{ \begin{array}{ll}
D_{t}^{\alpha}[x'(t)-g(t,x_{t})]=Ax(t)+f(t,x_{t})+\sigma(t,x_{t})\frac{\mathrm{d}w(t)}{\mathrm{d}t},\, t\in J,\, t\neq t_{i},\\
\triangle x(t_{i})=I_{i}(x_{t_{i}}),\,\triangle x'(t_{i})=J_{i}(x_{t_{i}}),\,i=1,2,\ldots,m,\\
x_{0}=\varphi\in\mathcal{B},\,x'_{0}=x_{1}\in H
\end{array} \right.
\end{eqnarray}
where $J=[0,b]$, $0<\alpha<1$, $D_{t}^{\alpha}$ denotes the Caputo fractional  derivative operator of order $\alpha$. $A:D(A)\subset H\rightarrow H$ is the infinitesimal generator of a sectorial operator. $g,\,f:J\times\mathcal{B}\rightarrow H$ and  $\sigma:J\times\mathcal{B}\rightarrow \mathcal{L}(G,H)$  are appropriate mappings. Here $\mathcal{B}$ is an abstract phase space to be defined later. The history $x_{t}:(-\infty,0]\rightarrow H,\, x_{t}(s)=x(t+s),s\leq0$ belongs to the abstract phase space $\mathcal{B}$.
Moreover, we denotes $\triangle x(t_{i})=x(t_{i}^{+})-x(t_{i}^{-})$ for $0\leq t_{0}<t_{1}<\cdots<t_{m}<t_{m+1}=b$ which are fixed numbers.
Let $x(t_{i}^{+})$ and $x(t_{i}^{-})$ represent the right and the left limits of $x(t)$ at $t=t_{i}$ respectively, Similarly $\triangle x'(t_{i})=x'(t_{i}^{+})-x'(t_{i}^{-})$ has the same meaning.

Let $(\Omega,\mathcal{F},P)$ be a complete probability space equipped with some filtration $\{\mathcal{F}_{t}\}_{t\geq0}$ satisfying the usual conditions, i.e., the filtration is right continuous and increasing while $\mathcal{F}_{0}$ contains all $P$-null sets. $H$, $G$ be two real separable Hilbert spaces. $<\cdot,\cdot>_{H}$,  $<\cdot,\cdot>_{G}$  denote the inner products on $H$ and $G$, respectively. And $|\cdot|_{H}$, $|\cdot|_{G}$  are vector norms on $H$, $G$. Let $\mathcal{L}(G,H)$ be the collection of all inner bounded operators from $G$ into $H$, with the usual operator norm $\|\cdot\|$.
The symbol $\{w(t),t\geq 0\}$ is a $G$ valued $\{\mathcal{F}_{t}\}_{t\geq0}$ Wiener process defined on the probability space $(\Omega,\mathcal{F},P)$ with covariance operator $Q$, i.e.
\begin{equation*}
    \mathbf{E}<w(t),x>_{G}<w(s),y>_{G}=(t\wedge s)<Qx,y>_{G},\:\forall x,y\in G
\end{equation*}
where $Q$ is a positive, self-adjoint and trace class operator on $K$. In particular, we regard $\{w(t),t\geq 0\}$ as a $G$ valued $Q$ wiener process related to $\{\mathcal{F}_{t}\}_{t\geq0}$ (see \cite{liuk, 16}), and $w(t)$ is defined as
\begin{equation*}
    w(t)=\sum_{n=1}^{\infty}\sqrt{\lambda_{n}}\beta_{n}(t)e_{n},\,t\geq0,
\end{equation*}
where $\beta_{n}(t)\,(n=1,2,3,\ldots)$ is a sequence of real valued standard Brownian  motions mutually independent on the probability space $(\Omega,\mathcal{F},P)$, let $\lambda_{n}$ $(n\in \mathbf{N})$ are the eigenvalues of $Q$ and $e_{n}$ $(n\in \mathbf{N})$ are the eigenvectors of  $\lambda_{n}$ corresponding to $\lambda_{n}$. That is
\begin{equation*}
    Qe_{n}=\lambda_{n}e_{n},\,n=1,2,3,\ldots.
\end{equation*}
In order to define stochastic integrals with respect to the $Q$ wiener process $w(t)$, we let $G_{0}=Q^{1/2}(G)$ of $G$ with the inner product,
\begin{equation*}
    <u,v>_{G_{0}}=<Q^{1/2}u,Q^{1/2}v>_{G}.
\end{equation*}
 It's easy to see $G_{0}$ is the subspace of $G$. Let $\mathcal{L}_{2}^{0}=\mathcal{L}_{2}(G_{0},H)$ denote the collection of all Hilbert Schmidt operators from $G_{0}$ into $H$. It turns out to be a separable Hilbert space equipped with the norm
\begin{equation*}
\|\psi\|^{2}_{\mathcal{L}_{2}^{0}}=tr\big((\psi Q^{1/2})(\psi Q^{1/2})^{*}\big), \quad\forall\,\psi\in \mathcal{L}_{2}^{0}.
\end{equation*}
Clearly, for any bounded operator $\psi\in \mathcal{L}(G,H)$, this norm reduces to $\|\psi\|^{2}_{\mathcal{L}_{2}^{0}}=tr(\psi Q\psi^{*})$.

Let $\Phi:(0,\infty)\rightarrow\mathcal{L}_{2}^{0}$ be a predictable and $\mathcal{F}_{t}$ adapted process such that
\begin{equation*}
    \int_{0}^{t}\mathbf{E}\|\Phi(s)\|^{2}_{\mathcal{L}_{2}^{0}}\mathrm{d}s<\infty,\, \forall t>0.
\end{equation*}
Then we can define the $H$ valued stochastic integral
\begin{equation*}
    \int_{0}^{t}\Phi(s)\mathrm{d}w(s),
\end{equation*}
which is a continuous square-integrable martingale \cite{0}. In the following we assume $\sigma:J\times\mathcal{B}\rightarrow \mathcal{L}_{2}^{0}$ in (\ref{eq1}). To the best of our knowledge this is the first time to consider the existence and uniqueness of the mild solution for Cauchy problem (\ref{eq1}).

The remainder of this paper is organized as follows. In section 2, we give some preliminaries which are used in this paper. In section 3, we give sufficient conditions for the existence and uniqueness of the mild solution of system (\ref{eq1}). In section 4, some sufficient conditions are introduced to guarantee the stability in mean square of the mild solution. In the last section, we present an example to support the results.

%

\section{Preliminaries}

The collection of all strongly measurable, square-integrable and $H$-valued random variables, denoted by $L_{2}(\Omega,H)$ is a Banach space equipped with norm $\|x(\cdot)\|_{L_{2}}=(\mathbf{E}\|x(\cdot)\|^{2}_{H})^{1/2}$ , where the expectation $\mathbf{E}$ is defined as $\mathbf{E}x=\int_{\Omega}x(\omega)\mathrm{d}P$.
Next, we present an axiomatic definition of the phase space $\mathcal{B}$ introduced in \cite{hale} and \cite{4}, where the axioms of the space $\mathcal{B}$ are established for $\mathcal{F}_{0}$-measurable functions from $(-\infty,0]$ into $H$, with a norm $\|\cdot\|_{\mathcal{B}}$ which satisfying
\begin{enumerate}
  \item [(A1)] If $x:(-\infty,b]\rightarrow H$, $b>0$ is such that $x_{0}\in\mathcal{B}$, then, for every $t\in J $, the following conditions hold:
      \begin{enumerate}
       \item [(1)] $x_{t}\in \mathcal{B},$
        \item [(2)] $|x(t)|\leq L\|x_{t}\|_{\mathcal{B}},$
        \item [(3)] $\|x_{t}\|_{\mathcal{B}}\leq \Gamma(t)\sup_{0\leq s\leq t}|x(s)|+N(t)\|x_{0}\|_{\mathcal{B}},$
        \end{enumerate}
        where $L>0$ is a constant; $\Gamma$, $N:[0,+\infty)\rightarrow[1,+\infty)$ are mappings. $\Gamma$ is continuous and $N$ is locally bounded. $L$, $\Gamma$, $N$ are independent on $x(\cdot)$.
  \item  [(A2)] The space $\mathcal{B}$ is complete.
\end{enumerate}
Then we have the following useful lemma (see \cite{4}).
\begin{lem}\label{lem}
Let $x:(-\infty,b]\rightarrow H$ be an $\mathcal{F}_{t}$ adapted measurable process such that the $\mathcal{F}_{0}$ adapted process $x_{0}=\varphi\in L_{2}(\Omega,\mathcal{B})$, then
\begin{equation}
    \mathbf{E}\|x_{s}\|_{\mathcal{B}}\leq N_{b}E\|\varphi\|_{\mathcal{B}}+\Gamma_{b}\mathbf{E}(\sup_{0\leq s\leq b}\|x(s)\|),
\end{equation}
where $N_{b}=\sup_{t\in J}\{N(t)\}$ and $\Gamma_{b}=\sup_{t\in J }\{\Gamma(t)\}.$
\end{lem}
\begin{defn}
Denote by $\mathcal{M}^{2}((-\infty,b],H)$ be the space of all $H$-valued c\`{a}dl\`{a}g measurable $\mathcal{F}_{t}$ adapted process $x=\{x(t)\}_{-\infty<t\leq b}$ such that
\begin{enumerate}
  \item [(i)] $x_{0}=\varphi\in\mathcal{B}$ and $x(t)$ is c\`{a}dl\`{a}g on $[0,b]$;
  \item [(ii)] endow the space $\mathcal{M}^{2}((-\infty,b],H)$ with the norm
  \begin{equation}\label{eqnorm}
    \|x\|_{\mathcal{M}^{2}}=\mathbf{E}\|\varphi\|_{\mathcal{B}}^{2}+\mathbf{E}(\sup_{t\in J}|x(t)|^{2})<\infty.
  \end{equation}
\end{enumerate}\end{defn}
 Then $\mathcal{M}^{2}((-\infty,b],H)$ with the norm (\ref{eqnorm}) is a Banach space, in the following of this paper, we use $\|\cdot\|$ for this norm.

\begin{defn}
A stochastic process $x(t):t\in(-\infty,b]\rightarrow H$ is called a mild solution of (\ref{eq1}) if
\begin{enumerate}
  \item [(i)] $x(t)$ is measurable and $\mathcal{F}_{t}$ adapted for all $t\in(-\infty,b],$ and $\{x_{t}:t\in[0,b]\}$ is $\mathcal{B}$-valued;
  \item [(ii)]$\int_{0}^{b}\|x(s)\|^{2}\mathrm{d}s<\infty$ ,$P$-a.s.;
  \item [(iii)]$x(t)$ has c\`{a}dl\`{a}g path on $t\in[0,b]$ a.s. and $x(t)$ satisfies the following integral equation for each $t\in[0,b]$,
  \begin{eqnarray}
   x(t) &=& S_{q}(t)\varphi(0)+\int_{0}^{t}S_{q} [x_{1}-g(0,\varphi)]\mathrm{d}s+\sum_{t_{i}<t}S_{q}(t-t_{i})I_{i}(x_{t_{i}})\\
   \nonumber &&+\sum_{t_{i}<t}\int_{t_{i}}^{t}S_{q}(t-s)[J_{i}(x_{t_{i}})-g(t_{i},x_{t_{i}}+I_{i}(x_{t_{i}}))+g(t_{i},x_{t_{i}})]\mathrm{d}s\\
    \nonumber &&+\int_{0}^{t}S_{q}(t-s)g(s,x_{s})\mathrm{d}s+\int_{0}^{t}T_{q}(t-s)f(s,x_{s})\mathrm{d}s\\
    \nonumber &&+\int_{0}^{t}T_{q}(t-s)\sigma(s,x_{s})\mathrm{d}w(s),
    \end{eqnarray}
    where $S_{q}(t)$, $T_{q}(t):\mathbf{R}_{+}\rightarrow \mathcal{L}(K,H)\, (q=1+\alpha)$ are given by
    \begin{eqnarray}
      S_{q}(t)&=& E_{q,1}(At^{q})=\frac{1}{2\pi i}\int_{B_{r}}\frac{e^{\lambda t}\lambda^{q-1}}{\lambda^{q}-A}\mathrm{d}\lambda, \\
     T_{q}(t) &=& t^{q-1}E_{q,q}(At^{q})=\frac{1}{2\pi i}\int_{B_{r}}\frac{e^{\lambda t}}{\lambda^{q}-A}\mathrm{d}\lambda,
    \end{eqnarray}
    and $B_{r}$ denotes the Bromwich path \cite{12};
   \item [(iv)] $x_{0}=\varphi\in\mathcal{B}.$
\end{enumerate}
\end{defn}
\begin{rem}
We should mention an important property of $S_{\gamma}(t)$ and $T_{\gamma}(t)$, that is there exist positive numbers $M$, $M_{b}$ such that $\|S_{\gamma}(t)\|_{\mathcal{L}(G,H)}\leq M$ and $\|T_{\gamma}(t)\|_{\mathcal{L}(G,H)}\leq t^{\gamma-1}M_{b}$ for $t\in J$, $\gamma\in(0,2)$ (\cite{3}), which plays an important role in the following discussion.
\end{rem}
\begin{lem}\label{B}(Bihari's inequality)
Assume $T>0$, $u_{0}\geq0$ and $u(t)$, $v(t)$ be continuous functions on $[0,T]$. Let $\kappa:\mathbf{R}_{+}\rightarrow\mathbf{R}_{+}$ be a concave continuous and nondecreasing function such that $\kappa(r)>0$ for all $r>0$. If
\begin{equation*}
    u(t)\leq u_{0}+\int_{0}^{t}v(s)\kappa(u(s))\mathrm{d}s \quad\textrm{for\,all}\,\quad 0\leq t\leq T,
\end{equation*}
then
\begin{equation*}
    u(t)\leq G^{-1}(G(u_{0})+\int_{0}^{t}v(s)\mathrm{d}s)
\end{equation*}
and for all $t\in [0,T]$, it holds that
\begin{equation*}
    G(u_{0})+\int_{0}^{t}v(s)\mathrm{d}s\in Dom(G^{-1}),
\end{equation*}
where $G(r)=\int_{1}^{r}\frac{\mathrm{d}s}{\kappa(s)},\, r\geq0 $ and $G^{-1}$ is the inverse function of $G$. In particular, if $u_{0}=0$ and  $\int_{0^{+}}\frac{\mathrm{d}s}{\kappa(s)}=\infty$, then $u(t)=0$ for all $0\leq t\leq T$.
\end{lem}
 The following lemma is useful in the proof of the exponential stability of the mild solution which is an analogue of Theorem 18 of \cite{4}.
\begin{lem}\label{cc}\cite{bh}
Let the assumption of Lemma \ref{B} hold and $v(t)\geq0$ for all $t\in[0,T].$ If for all $\epsilon>0$, there exists $t_{1}\geq0$ for all $0\leq u_{0}\leq\epsilon,$ $\int_{t_{1}}^{T}v(s)\mathrm{d}s\leq\int_{u_{0}}^{\epsilon}\frac{1}{\kappa(s)}\mathrm{d}s$ holds. Then for every $t\in[t_{1},T],$ the estimates $u(t)\leq \epsilon$ holds.
\end{lem}

\section{Existence of the mild solution }
 In this section we first make the following hypotheses.
 \begin{enumerate}
   \item [(H1)] $g,f:\,J\times\mathcal{B}\rightarrow H$ and $\sigma:J\times\mathcal{B}\rightarrow H$ satisfy \\
        $|g(t,\varphi)-g(t,\phi)|^{2}\vee|f(t,\varphi)-f(t,\phi)|^{2}\vee|\sigma(t,\varphi)-\sigma(t,\phi)|^{2}\leq \kappa(\|\varphi-\phi\ |^{2}_{\mathcal{B}}),$ for all $t\in J$,$\varphi,\phi\in\mathcal{B}$,
        where $\kappa(\cdot)$ is a concave, nondecreasing and continuous function from $\mathbf{R}_{+}\rightarrow\mathbf{R}_{+}$ such that $\kappa(0)=0$, $\kappa(u)>0$ for $u>0$ and $\int_{0^{+}}\frac{\mathrm{d}s}{\kappa(s)}=\infty.$
   \item [(H2)] $I_{k},J_{k}: \mathcal{B}\rightarrow H$ are continuous and there are positive constants $p_{k}, \,q_{k}(k=1,2,\ldots,m) $ such that for each $\varphi,\phi\in\mathcal{B},$
\begin{equation*}
    |I_{k}(\varphi)-I_{k}(\phi)|^{2}\leq p_{k}\|\varphi-\phi\|^{2},\, |J_{k}(\varphi)-J_{k}(\phi)|^{2}\leq q_{k}\|\varphi-\phi\|^{2}.
\end{equation*}
\item [(H3)]$|g(t,0)|^{2}\vee|f(t,0)|^{2}\vee|\sigma(t,0)|^{2}\leq K,$ $K$ is a positive constant, and $I_{k}(0)=0,\,J_{k}(0)=0,k=1,2,\ldots,m.$
 \end{enumerate}
We consider the sequence of successive approximations defined as follows:
\begin{equation}
x^{0}(t)=S_{q}(t)\varphi(0)+\int_{0}^{t}S_{q}(s)[x_{1}-g(0,\varphi)]\mathrm{d}s,\,t\in J,
\end{equation}
\begin{eqnarray}\label{eqmain}
\nonumber &&x^{n}(t)=S_{q}(t)\varphi(0)+\int_{0}^{t}S_{q}(s)[x_{1}-g(0,\varphi)]\mathrm{d}s+\sum_{t_{i}<t}S_{q}(t-t_{i})I_{i}(x_{t_{i}}^{n-1})\\
  \nonumber&&\quad+\sum_{t_{i}<t}\int_{t_{i}}^{t}S_{q}(t-s)[J_{i}(x_{t_{i}}^{n-1})-g(t_{i},x_{t_{i}}^{n-1}+I_{i}(x^{n-1}_{t_{i}}))+g(t_{i},x_{t_{i}}^{n-1})]\mathrm{d}s\\
  \nonumber&&\quad+\int_{0}^{t}S_{q}(t-s)g(s,x_{s}^{n-1})\mathrm{d}s+\int_{0}^{t}T_{q}(t-s)f(s,x_{s}^{n-1})\mathrm{d}s\\
  &&\quad+\int_{0}^{t}T_{q}(t-s)\sigma(s,x_{s}^{n-1})\mathrm{d}w(s),\,t\in J, n\geq1,
\end{eqnarray}
\begin{equation}
 x^{n}(t)=\varphi(t),\,-\infty<t\leq0, n\geq1.
 \end{equation}
\begin{lem}\label{lembb}
Assume the (H1)-(H3) hold, and
\begin{equation*}
7mM^{2}\Gamma_{b}\sum_{i=1}^{m}p_{i}+14mM^{2}b^{2}\Gamma_{b}\sum_{i=1}^{m}q_{i}<1,
\end{equation*}
  then $x^{n}(t)\in \mathcal{M}^{2}((-\infty,b);H)$ for all $t\in(-\infty,b]$, $n\geq0$, that is
\begin{equation}
    \mathbf{E}\|x^{n}(t)\|^{2}\leq \tilde{M},\,n=1,2\ldots.
\end{equation}
where $\tilde{M}$ is a positive constant.
\end{lem}
\begin{proof} Obviously $x^{0}(t)\in \mathcal{M}^{2}((-\infty,b),H)$, and
\begin{eqnarray*}
  &&\mathbf{E}|X^{n}(t)|^{2} \\
  &&\leq 7\mathbf{E}|S_{q}(t)\varphi(0)|^{2}+ 7\mathbf{E}|\int_{0}^{t}S_{q}(s)[x_{1}-g(0,\varphi)]\mathrm{d}s|^{2}\\
  &&\quad+7\mathbf{E}|\sum_{t_{i}<t}S_{q}(t-t_{i})I_{i}(x_{t_{i}}^{n-1})|^{2}\\
   && \quad+7\mathbf{E}|\sum_{t_{i}<t}\int_{t_{i}}^{t}S_{q}(t-s)[J_{i}(x_{t_{i}}^{n-1})-g(t_{i},x_{t_{i}}^{n-1}+I_{i}(x^{n-1}_{t_{i}}))+g(t_{i},x_{t_{i}}^{n-1})]\mathrm{d}s|^{2}\\
     &&\quad+7\mathbf{E}|\int_{0}^{t}S_{q}(t-s)g(s,x_{s}^{n-1})\mathrm{d}s|^{2}+7E|\int_{0}^{t}T_{q}(t-s)f(s,x_{s}^{n-1})\mathrm{d}s)|^{2}\\
    && \quad+7\mathbf{E}|\int_{0}^{t}T_{q}(t-s)\sigma(s,x_{s}^{n-1})\mathrm{d}w(s)|^{2}\\
     &&=\Lambda_{1}+\Lambda_{2}+\Lambda_{3}+\Lambda_{4}+\Lambda_{5}+\Lambda_{6}+\Lambda_{7}.
\end{eqnarray*}
It's easy to get the estimations
$\Lambda_{1}\leq 7M^{2}\mathbf{E}|\varphi(0)|^{2}$, $\Lambda_{2}\leq 21 M^{2}b^{2}(|x_{1}|^{2}+\kappa(\|\varphi\|_{\mathcal{B}}^{2})+K)$, and
\begin{equation*}
   \Lambda_{3}\leq 7mM^{2}\sum_{t_{i}<t}\mathbf{E}\|I_{i}(x_{t_{i}}^{n-1})\|_{\mathcal{B}}^{2}\leq 7mM^{2}\sum_{t_{i}<t}p_{i}\mathbf{E}\|x_{t_{i}}^{n-1}\|_{\mathcal{B}}.
\end{equation*}
By the fact $\|S_{\gamma}(t)\|_{\mathcal{L}(G,H)}\leq M$ and (H2)-(H3), we have
\begin{eqnarray*}
  \Lambda_{4}
  &\leq&14\mathbf{E}|\sum_{t_{i}<t}\int_{t_{i}}^{t}S_{q}(t-s)J_{i}(x_{t_{i}}^{n})\mathrm{d}s |^{2}\\
  &&+14\mathbf{E}|\sum_{t_{i}<t}\int_{t_{i}}^{t}S_{q}(t-s)(g(t_{i},x_{t_{i}}^{n-1}+I_{i}(x_{t_{i}}^{n-1}))-g(t_{i},x_{t_{i}}^{n-1}))\mathrm{d}s |^{2}\\
   &\leq&14mM^{2}b\sum_{t_{i}<t}\mathbf{E}\int_{t_{i}}^{t}|J_{i}(x_{t_{i}})|^{2}\mathrm{d}s+14mM^{2}b\sum_{t_{i}<t}\mathbf{E}\int_{t_{i}}^{t}\kappa(|I(x_{t_{i}}^{n-1})|^{2})\mathrm{d}s \\ &\leq&14mM^{2}b^{2}\sum_{t_{i}<t}q_{i}\mathbf{E}\|x_{t_{i}}^{n-1}\|_{\mathcal{B}}^{2}+14mM^{2}b\sum_{t_{i}<t}\int_{t_{i}}^{t}\kappa(\mathbf{E}(p_{i}\|x_{t_{i}}^{n-1})\|_{\mathcal{B}}^{2})\mathrm{d}s. \end{eqnarray*}
   and
\begin{eqnarray*}
 \Lambda_{5}&\leq& 7\mathbf{E}|\int_{0}^{t}S_{q}(t-s)g(s,x_{s}^{n-1})\mathrm{d}s|^{2}\\
 &\leq&7M^{2}b \mathbf{E}\int_{0}^{t}|g(s,x_{s}^{n-1})-g(s,0)+g(s,0)|^{2}\mathrm{d}s\\
&\leq&14M^{2}b\mathbf{E}\int_{0}^{t}[|g(s,x_{s}^{n-1})-g(s,0)|^{2}+|g(s,0)|^{2}]\mathrm{d}s\\
&\leq&14M^{2}b\int_{0}^{t}\kappa(\|x_{s}^{n-1}\|_{\mathcal{B}}^{2})\mathrm{d}s+14M^{2}b^{2}K.
  \end{eqnarray*}
  Since $\|T_{\gamma}(t)\|_{\mathcal{L}(G,H)}\leq t^{\gamma-1}M_{b}$ and (H1)-(H3), we get the following inequality
\begin{eqnarray*}
\Lambda_{6}&\leq& 7\mathbf{E}|\int_{0}^{t}T_{q}(t-s)f(s,x_{s}^{n-1})\mathrm{d}s|^{2}\\
&\leq&7M_{b}^{2}\frac{b^{2q-1}}{2q-1}\mathbf{E}\int_{0}^{t}|f(s,x_{s}^{n-1})-f(s,0)+f(s,0)|^{2}\mathrm{d}s\\
&\leq&14M_{b}^{2}\frac{b^{2q-1}}{2q-1}\int_{0}^{t}\kappa(\mathbf{E}\|x_{s}^{n-1}\|_{\mathcal{B}}^{2})\mathrm{d}s+14M_{b}^{2}\frac{b^{2q}}{2q-1}K.
  \end{eqnarray*}
  We apply the H\"{o}lder inequality and the Burkholder-Davis-Gundy inequality to $\Lambda_{7}$, combining(H1)-(H3), we can obtain
\begin{eqnarray*}
\Lambda_{7}&\leq& 7\mathbf{E}|\int_{0}^{t}T_{q}(t-s)\sigma(s,x_{s}^{n-1})\mathrm{d}w(s)|^{2}\\
&\leq&7\mathbf{E}\int_{0}^{t}|T_{q}(t-s)\sigma(s,x_{s}^{n-1})|^{2}\mathrm{d}s\\
&\leq&7M_{b}^{2}b^{2q-2}\mathbf{E}\int_{0}^{t}|\sigma(s,x_{s}^{n-1})-\sigma(s,0)+\sigma(s,0)|^{2}\mathrm{d}s\\
&\leq&14M_{b}^{2}b^{2q-2}\int_{0}^{t}\kappa(\mathbf{E}\|x_{s}^{n-1}\|_{\mathcal{B}}^{2})\mathrm{d}s+14M_{b}^{2}b^{2q-1}K.
\end{eqnarray*}
    Let\begin{eqnarray*}
        c_{1} &=& 7M^{2}\mathbf{E}|\varphi(0)|^{2}+21 M^{2}b^{2}(|x_{1}|^{2}+\kappa(\|\varphi\|_{\mathcal{B}}^{2})+K) \\
         &&  +14M^{2}b^{2}K+14M_{b}^{2}\frac{b^{2q}}{2q-1}K+14M_{b}^{2}b^{2q-1},
       \end{eqnarray*}
the estimations for $\Lambda_{i}$ (i=1,2,\ldots,7) together yields
\begin{eqnarray*}
  \mathbf{E}|X^{n}(t)|^{2} &\leq& c_{1}+7mM^{2}\sum_{t_{i}<t}p_{i}\mathbf{E}\|x_{t_{i}}^{n-1}\|_{\mathcal{B}}+ 14mM^{2}b^{2}\sum_{t_{i}<t}q_{i}\mathbf{E}\|x_{t_{i}}^{n-1}\|_{\mathcal{B}}^{2}\\&&+14mM^{2}b\sum_{t_{i}<t}\int_{t_{i}}^{t}\kappa(\mathbf{E}(p_{i}\|x_{t_{i}}^{n-1})\|_{\mathcal{B}}^{2})\mathrm{d}s\\
  &\,&+(14M^{2}b+14M_{b}^{2}\frac{b^{2q-1}}{2q-1}+14M_{b}^{2}b^{2q-2})\int_{0}^{t}\kappa(\|x_{s}^{n-1}\|_{\mathcal{B}}^{2})\mathrm{d}s.
\end{eqnarray*}
By Lemma \ref{lem} and the property of $\kappa(\cdot)$, we can find a pair of positive constants $\alpha$ and $\beta$, such that $\kappa(u)\leq \alpha+\beta u,\quad\forall u\geq0.$ Then
\begin{eqnarray*}
 && \mathbf{E}\sup_{0\leq s\leq t}|X^{n}(s)|^{2}\\
   &\leq& c_{1}+7mM^{2}N_{b}\sum_{t_{i}<t}p_{i}\mathbf{E}\|\varphi\|_{\mathcal{B}}+14mM^{2}b^{2}\sum_{t_{i}<t}p_{i}\mathbf{E}\|\varphi\|_{\mathcal{B}}+14m^{2}M^{2}b^{2}\alpha\\
   &\,&+(14M^{2}b+14M_{b}^{2}\frac{b^{2q-1}}{2q-1}+14M_{b}^{2}b^{2q-2})b\alpha\\
&\,&+(7mM^{2}\Gamma_{b}\sum_{t_{i}<t}p_{i}+14mM^{2}b^{2}\Gamma_{b}\sum_{t_{i}<t}q_{i})\mathbf{E}\sup_{0\leq s\leq t}|X^{n-1}(s)|^{2}\\&&+(14mM^{2}b \sum_{t_{i}<t}p_{i}+14M^{2}b+14M_{b}^{2}\frac{b^{2q-1}}{2q-1}+14M_{b}^{2}b^{2q-2})\\
 &\,&\quad\times\beta\mathbf{E}\int_{0}^{t}\sup_{0\leq\theta\leq s}|x^{n-1}(\theta)|\mathrm{d}s,
\end{eqnarray*}
and
\begin{eqnarray*}
&&\max_{1\leq n\leq k}\{\mathbf{E}\sup_{0\leq s\leq t}|x^{n}(s)|^{2}\} \\
&&\leq c_{1}+7mM^{2}N_{b}\sum_{t_{i}<t}p_{i}\mathbf{E}\|\varphi\|_{\mathcal{B}}+14mM^{2}b^{2}\sum_{t_{i}<t}p_{i}\mathbf{E}\|\varphi\|_{\mathcal{B}}+14m^{2}M^{2}b\alpha\\
&&\quad+(14M^{2}b+14M_{b}^{2}\frac{b^{2q-1}}{2q-1}+14M_{b}^{2}\frac{b^{2q-1}}{2q-1})b\alpha\\
&&\quad+(7mM^{2}\Gamma_{b}\sum_{t_{i}<t}p_{i}+14mM^{2}b^{2}\Gamma_{b}\sum_{t_{i}<t}q_{i})\max_{1\leq n\leq k}\mathbf{E}\sup_{0\leq s\leq t}|X^{n}(s)|^{2}\\
&&\quad+(14mM^{2}b \sum_{t_{i}<t}p_{i}+14M^{2}b+14M_{b}^{2}\frac{b^{2q-1}}{2q-1}+14M_{b}^{2}\frac{b^{2q-1}}{2q-1})\\
&&\quad\quad\times\beta\int_{0}^{t}\max_{1\leq n\leq k}\mathbf{E}\sup_{0\leq s\leq t}|x^{n}(s)|^{2}\mathrm{d}s.
\end{eqnarray*}
If we let
\begin{eqnarray*}
c_{2}
&=&\frac{c_{1}+7mM^{2}N_{b}\sum_{t_{i}<t}p_{i}\mathbf{E}\|\varphi\|_{\mathcal{B}}+14mM^{2}b^{2}\sum_{t_{i}<t}p_{i}\mathbf{E}\|\varphi\|_{\mathcal{B}}}{1-7mM^{2}\Gamma_{b}\sum_{t_{i}<t}p_{i}-14mM^{2}b^{2}\Gamma_{b}\sum_{t_{i}<t}q_{i}}\\
&\quad&+\frac{14m^{2}M^{2}b\alpha+(14M^{2}b^{2}+14M_{b}^{2}b^{2q}+14M_{b}^{2}b^{2q})\alpha}{1-7mM^{2}\Gamma_{b}\sum_{t_{i}<t}p_{i}-14mM^{2}b^{2}\Gamma_{b}\sum_{t_{i}<t}q_{i}}
\end{eqnarray*}
\begin{equation*}
c_{3}=\frac{(14 m M^{2}b \sum_{t_{i}<t}p_{i}+14M^{2}b+14M_{b}^{2}\frac{b^{2q-1}}{2q-1}+14M_{b}^{2}b^{2q-2})\beta}{1-7mM^{2}\Gamma_{b}\sum_{t_{i}<t}p_{i}-14mM^{2}b^{2}\Gamma_{b}\sum_{t_{i}<t}q_{i}};
\end{equation*}
then
\begin{eqnarray}
  \max_{1\leq n\leq k}\{\mathbf{E}\sup_{0\leq s\leq t}|X^{n}(s)|^{2}\}\leq c_{2}+c_{3} \int_{0}^{t}\max_{1\leq n\leq k}\mathbf{E}\sup_{0\leq s\leq t}|x^{n}(s)|^{2}\mathrm{d}s,
\end{eqnarray}
 by the Gronwall inequality we have
\begin{equation*}
   \max_{1\leq n\leq k}\{\mathbf{E}\sup_{0\leq s\leq t}|x^{n}(s)|^{2}\}\leq c_{2}e^{c_{3}}.
\end{equation*}
 Due to the arbitrary of $k$, we have
\begin{equation*}
\mathbf{E}\sup_{0\leq s\leq t}|x^{n}(s)|^{2}\leq c_{2}e^{c_{3}},\, \textrm{for\,all}\,0\leq t\leq b,\,n\geq1.
\end{equation*}
Consequently,
\begin{equation}
    \|x^{n}(t)\|^{2}\leq\mathbf{ E}\|\varphi\|_{\mathcal{B}}^{2}+\mathbf{E}(\sup_{0\leq s\leq b}|x^{n}(s)|^{2})\leq \mathbf{E}\|\varphi\|_{\mathcal{B}}^{2}+ c_{2}e^{c_{3}}<\infty,
\end{equation}
so we can take $\tilde{M}=\mathbf{E}\|\varphi\|_{\mathcal{B}}^{2}+ c_{2}e^{c_{3}}$. This completes the proof of Lemma \ref{lembb}.\end{proof}
\begin{thm}\label{theorem}
If (H1)-(H3) and
\begin{equation}\label{1}
    \max\{7mM^{2}\Gamma_{b}\sum_{i=1}^{m}p_{i}+14mM^{2}b^{2}\Gamma_{b}\sum_{i=1}^{m}q_{i},7mM^{2}\sum_{i=1}^{m}p_{i}+7mM^{2}b\sum_{i=1}^{m}q_{i}\}<1
\end{equation}
 hold, then the Cauchy problem (\ref{eq1}) has a unique mild solution on $(-\infty,b]$.
\end{thm}
\begin{proof} Since
\begin{eqnarray*}
 &&|x^{n+l}(t)-x^{n}(t)|^{2}\\
&&= |\sum_{t_{i}<t}S_{q}(t-t_{i})[I_{i}(x_{t_{i}}^{l+n-1})-I_{i}(x_{t_{i}}^{n-1})] \\
 &&\quad +\sum_{t_{i}<t}\int_{t_{i}}^{t}S_{q}(t-s)[J(x_{t_{i}}^{l+n-1})-J(x_{t_{i}}^{n-1})]\mathrm{d}s\\
 &&\quad-\sum_{t_{i}<t}\int_{t_{i}}^{t}S_{q}(t-s)[ g(t_{i},x_{t_{i}}^{l+n-1}+I_{i}(x_{t_{i}}^{l+n-1}))-g(t_{i},x_{t_{i}}^{n-1}+I_{i}(x_{t_{i}}^{n-1}))\\
 &&\quad\quad+g(t_{i},x_{t_{i}}^{l+n-1})-g(t_{i},x_{t_{i}}^{n-1})]\mathrm{d}s\\
   &&\quad+\int_{0}^{t}S_{q}(t-s)[g(s,x_{s}^{n+l-1})-g(s,x_{s}^{n-1})]+T_{q}(t-s)[f(s,x_{s}^{n+l-1})-f(s,x_{s}^{n-1})]\mathrm{d}s\\
   &&\quad+\int_{0}^{t}T_{q}(t-s)[\sigma(s,x_{s}^{n+l-1})-\sigma(s,x_{s}^{n-1})]\mathrm{d}w(s)|^{2}.
\end{eqnarray*}
By the fact $\|S_{\gamma}(t)\|_{\mathcal{L}(G,H)}\leq M$, $\|T_{\gamma}(t)\|_{\mathcal{L}(G,H)}\leq t^{\gamma-1}M_{b}$ for $t\in J$ and (H1)-(H3), we get
\begin{eqnarray*}
&&|x^{n+l}(t)-x^{n}(t)|^{2}\\
&&\leq7|\sum_{t_{i}<t}S_{q}(t-t_{i})[I_{i}(x_{t_{i}}^{n+l-1})-I_{i}(x_{t_{i}}^{n-1})]|^{2}\\
&&+7|\sum_{t_{i}<t}\int_{t_{i}}^{t}S_{q}(t-s)[J(x_{t_{i}}^{n+l-1})-J(x_{t_{i}}^{n-1})]\mathrm{d}s|^{2}\\
&&+7|\sum_{t_{i}<t}\int_{t_{i}}^{t}S_{q}(t-s)[ g(t_{i},x_{t_{i}}^{l+n-1}+I_{i}(x_{t_{i}}^{l+n-1}))-g(t_{i},x_{t_{i}}^{n-1}+I_{i}(x_{t_{i}}^{n-1}))]\mathrm{d}s|^{2}
\\
&&+7|\sum_{t_{i}<t}\int_{t_{i}}^{t}S_{q}(t-s)[g(t_{i},x_{t_{i}}^{l+n-1})-g(t_{i},x_{t_{i}}^{n-1})]\mathrm{d}s|^{2}\\
&&+7|\int_{0}^{t}S_{q}(t-s)[g(s,x_{s}^{n+l-1})-g(s,x_{s}^{n-1})]\mathrm{d}s|^{2}\\&&+|\int_{0}^{t}T_{q}(t-s)[f(s,x_{s}^{n+l-1})-f(s,x_{s}^{n-1})]\mathrm{d}s|^{2}\\
&&+7|\int_{0}^{t}T_{q}(t-s)[\sigma(s,x_{s}^{n+l-1})-\sigma(s,x_{s}^{n-1})]\mathrm{d}w(s)|^{2}
\end{eqnarray*}
\begin{eqnarray*}
&\leq&7mM^{2}\sum_{t_{i}<t}p_{i}\sup_{0\leq s\leq t}|x^{l+n-1}(s)-x^{n-1}(s)|^{2}\\
&&+7mM^{2}b\sum_{t_{i}<t}q_{i}\sup_{0\leq s\leq t}|x^{l+n-1}(s)-x^{n-1}(s)|^{2}\\
&&+7mM^{2}\sum_{t_{i}<t}\int_{t_{i}}^{t}\kappa((1+p_{i})\sup_{0\leq s\leq t_{i}}|x^{l+n-1}(s)-x^{n-1}(s)|^{2})\mathrm{d}s\\
&&+7mM^{2}\sum_{t_{i}<t}\int_{t_{i}}^{t}\kappa(\sup_{0\leq s\leq t_{i}}|x^{l+n-1}(s)-x^{n-1}(s)|^{2})\mathrm{d}s\\
&&+7M^{2}b\int_{0}^{t}\kappa(\sup_{0\leq r\leq s}|x^{l+n-1}(r)-x^{n-1}(r)|^{2})\mathrm{d}s\\
&&+7M_{b}^{2}\frac{b^{2q-1}}{2q-1}\int_{0}^{t}\kappa(\sup_{0\leq r\leq s}|x^{l+n-1}(r)-x^{n-1}(r)|^{2})\mathrm{d}s\\
&&+7M_{b}^{2}b^{2q-2}\int_{0}^{t}\kappa(\sup_{0\leq r\leq s}|x^{l+n-1}(r)-x^{n-1}(r)|^{2})\mathrm{d}s\\
&=&(7mM^{2}\sum_{t_{i}<t}p_{i}+7mM^{2}b\sum_{t_{i}<t}q_{i})\sup_{0\leq s\leq t}|x^{l+n-1}(s)-x^{n-1}(s)|^{2}\\
&&+7mM^{2}\sum_{t_{i}<t}\int_{t_{i}}^{t}\kappa((1+p_{i})\sup_{0\leq r\leq t_{i}}|x^{l+n-1}(r)-x^{n-1}(r)|^{2})\mathrm{d}s\\
&&+(7m^{2}M^{2}+7M^{2}b+7M_{b}^{2}\frac{b^{2q-1}}{2q-1}+7M_{b}^{2}b^{2q})\int_{0}^{t}\kappa(\sup_{0\leq r\leq s}|x^{l+n-1}(r)-x^{n-1}(r)|^{2})\mathrm{d}s.
\end{eqnarray*}
Let $\bar{p}=\max_{1\leq i\leq m}{p_{i}}$. Since
\begin{eqnarray}
  &&\nonumber\int_{t_{i}}^{t}\kappa\big((1+p_{i})\sup_{0\leq r\leq t_{i}}|x^{l+n-1}(r)-x^{n-1}(r)|^{2}\big)\mathrm{d}s\\
  \nonumber&&\leq  \int_{t_{i}}^{t}\kappa\big((1+p_{i})\sup_{0\leq r\leq s}|x^{l+n-1}(r)-x^{n-1}(r)|^{2}\big)\mathrm{d}s \\
  &&\leq \int_{0}^{t}\kappa\big((1+\bar{p})\sup_{0\leq r\leq s}|x^{l+n-1}(r)-x^{n-1}(r)|^{2}\big)\mathrm{d}s .
\end{eqnarray}
Obviously $\tilde{\kappa}\circ a(\cdot)=\kappa(a(\cdot))$ is also a concave function, then we get
\begin{eqnarray}\label{kk}
\nonumber&&\sup_{0\leq s\leq t}|x^{n+l}(t)-x^{n}(t)|^{2}\\&&\leq(
7mM^{2}\sum_{t_{i}<t}p_{i}+7mM^{2}b\sum_{t_{i}<t}q_{i})\sup_{0\leq s\leq t}|x^{n+l-1}(s)-x^{n-1}(s)|^{2}\\
\nonumber&&+(14m^{2}M^{2}+7M^{2}b+7M_{b}^{2}\frac{b^{2q-1}}{2q-1}+7M_{b}^{2}b^{2q-2})\\
\nonumber&&\times\int_{0}^{t}\tilde{\kappa}(\sup_{0\leq r\leq s}|x^{l+n-1}(r)-x^{n-1}(r)|^{2})\mathrm{d}s.
\end{eqnarray}
It is easy to see
\begin{eqnarray}
\nonumber&&|x^{n+l}(t)-x^{n}(t)|^{2}-(7mM^{2}\sum_{t_{i}<t}p_{i}+7mM^{2}M^{2}b\sum_{t_{i}<t}q_{i})|x^{l+n-1}(r)-x^{n-1}(r)|^{2}\\
 && \leq(14m^{2}M^{2}+7M^{2}b+7M_{b}^{2}\frac{b^{2q-1}}{2q-1}+7M_{b}^{2}b^{2q-2})\\
&&\nonumber \quad\times\int_{0}^{t}\tilde{\kappa}(\sup_{0\leq r\leq s}|x^{l+n-1}(r)-x^{n-1}(r)|^{2})\mathrm{d}s.
\end{eqnarray}
From lemma \ref{lembb}, we get
\begin{eqnarray}
\nonumber &&|x^{n+l}(t)-x^{n}(t)|^{2}-(7mM^{2}\sum_{t_{i}<t}p_{i}+7mM^{2}b\sum_{t_{i}<t}q_{i})|x^{l+n-1}(r)-x^{n-1}(r)|^{2} \\
 \nonumber &&\leq(14m^{2}M^{2}+7M^{2}b+7M_{b}^{2}\frac{b^{2q-1}}{2q-1}+7M_{b}^{2}b^{2q-2})\int_{0}^{t}\tilde{\kappa}(2M_{1})\mathrm{d}s\\
  &&\leq c_{4}\tilde{\kappa}(2M_{1})t=c_{5}t.
\end{eqnarray}
Define
\begin{eqnarray}
  \varphi_{1}(t) =c_{5}t,\, \varphi_{n+1}(t)=c_{4} \int_{0}^{t}\tilde{\kappa}(\varphi_{n}(s))\mathrm{d}s,\,n\geq1.
\end{eqnarray}
Choose $b_{1}\in [0,b)$ such that $c_{4}\tilde{\kappa}(c_{5}t)\leq c_{5}$, for all $0\leq t\leq b_{1}$.

 We give the statement that for any $t\in [0,b_{1})$, $\{\varphi_{n}(t)\}$ is a decreasing sequence. In fact
\begin{eqnarray*}
 \varphi_{2}(t)=c_{4}\int_{0}^{t}\tilde{\kappa}(\varphi_{1}(s))\mathrm{d}s=c_{4}\int_{0}^{t}\tilde{\kappa}(c_{5}s)\mathrm{d}s \leq \int_{0}^{t} c_{5}\mathrm{d}s=\varphi_{1}(t).
\end{eqnarray*}
 By induction, we get
 \begin{eqnarray}
 \varphi_{n+1}(t)=c_{4}\int_{0}^{t}\tilde{\kappa}(\varphi_{n}(s))\mathrm{d}s\leq c_{4}\int_{0}^{t}\tilde{\kappa}(\varphi_{n-1}(s))\mathrm{d}s=\varphi_{n}(t)\quad \forall \,0\leq t\leq b_{1}.
 \end{eqnarray}
 Therefore, the statement is true, and we can define the function $\phi(t)$ as
 \begin{eqnarray}
 \phi(t)&=&\lim_{n\rightarrow\infty}\varphi_{n}(t)=\lim_{n\rightarrow\infty}c_{4} \int_{0}^{t}\tilde{\kappa}(\varphi_{n-1}(s))\mathrm{d}s\\\nonumber&=&\lim_{n\rightarrow\infty}c_{4} \int_{0}^{t}\tilde{\kappa}(\phi(s))\mathrm{d}s,\,0\leq t\leq b_{1};
 \end{eqnarray}
 By the Bihari's inequality, we get $ \phi(t)=0$ for all  $0\leq t\leq b_{1}$. It means that for all $0\leq t\leq b_{1}$
\begin{eqnarray}\label{eqlim}
&&\nonumber\lim_{n\rightarrow\infty}[\mathbf{E}|x^{n+l}(t)-x^{n}(t)|^{2}\\
&&\quad-(7mM^{2}\sum_{t_{i}<t}p_{i}+7mM^{2}b\sum_{t_{i}<t}q_{i})\mathbf{E}|x^{m+l-1}(t)-x^{n-1}(t)|^{2}]=0.
\end{eqnarray}
Using the assumption of the theorem $7mM^{2}\sum_{t_{i}<t}p_{i}+7mM^{2}b\sum_{t_{i}<t}q_{i}<1$ and (\ref{eqlim}), we get
\begin{eqnarray}
\lim_{n\rightarrow\infty}\mathbf{E}|x^{n+l}(t)-x^{n}(t)|^{2}=0,\,0\leq t\leq b_{1},
\end{eqnarray}
which means that $\{x^{n}(t)\}$ is a Cauchy sequence in $L^{2}$. Let
\begin{equation*}
 \lim_{n\rightarrow\infty}x^{n}(t)=x(t),
\end{equation*}
 obviously
\begin{equation*}
 \|x(t)\|^{2}\leq \tilde{M},\,0\leq t\leq b_{1}.
\end{equation*}
Taking limits on both side of equation (\ref{eqmain}), for all $t\in [0,b_{1}]$, we have
\begin{eqnarray}\label{eqmain1}
 x(t)&=&S_{q}(t)\varphi(0)+\int_{0}^{t}S_{q}(s)[x_{1}-g(0,\varphi)]\mathrm{d}s+\sum_{t_{i}<t}S_{q}(t-t_{i})I_{i}(x_{t_{i}})\\
\nonumber&&+\sum_{t_{i}<t}\int_{t_{i}}^{t}S_{q}(t-s)[J_{i}(x_{t_{i}})-g(t_{i},x_{t_{i}}+I_{i}(x_{t_{i}}))+g(t_{i},x_{t_{i}})]\mathrm{d}s\\
 \nonumber &&+\int_{0}^{t}S_{q}(t-s)g(s,x_{s})\mathrm{d}s+\int_{0}^{t}T_{q}(t-s)f(s,x_{s})\mathrm{d}s\\
 \nonumber &&
 +\int_{0}^{t}T_{q}(t-s)\sigma(s,x_{s})\mathrm{d}w(s).
\end{eqnarray}
So we have presented the existence of the mild solution of problem (\ref{eq1}) on $[0,b_{1}] $. By iteration we can get the existence of the mild solution of  problem (\ref{eq1}) on $[0,b]$.

Suppose that $x(t)$ and $\bar{x}(t)$ are two solutions of (\ref{eq1}). In the similar discussion as (\ref{eqlim}) we can get
\begin{eqnarray}
 &&[1-(7mM^{2}\sum_{t_{i}<t}p_{i}+7mM^{2}b\sum_{t_{i}<t}q_{i})]\mathbf{E}|x(t)-\bar{x}(t)|^{2} \\\nonumber&&\leq
 (14m^{2}M^{2}+7M^{2}b+7M_{b}^{2}\frac{b^{2q-1}}{2q-1}+7M_{b}^{2}b^{2q-2})\\\nonumber&&\quad\times\int_{0}^{t}\tilde{\kappa}(\mathbf{E}\sup_{0\leq r\leq s}|x^{m+n-1}(r)-x^{n-1}(r)|^{2})\mathrm{d}s,
\end{eqnarray}
the  Bihari inequality implies $\mathbf{E}|x(t)-\bar{x}(t)|^{2}=0$, and we have show the existence and uniqueness of the mild solution of (\ref{eq1}).\end{proof}

\section{stability of solutions}
In this section, we will give the continuous dependence of solutions on the initial value by means of the Bihari's inequality. We first propose the following the assumption on $g$ instead of (H1),
\begin{enumerate}
  \item [(H4)] $g:\,J\times\mathcal{B}\rightarrow H$ satisfies  $|g(t,\varphi)-g(t,\phi)|^{2}\leq K_{1}\|\varphi-\phi\|_{\mathcal{B}}^{2}$, where $K_{1}$ is a positive constant.
\end{enumerate}
\begin{defn}\cite{4}
Assume a mild solution $x^{\varphi,x_{1}}(t)$ of Cauchy problem (\ref{eq1}) with initial value $(\varphi,x_{1})$ stable is said to be stable in square if for all $\epsilon>0$ there exists $\delta>0$ such that
\begin{equation}
  \mathbf{E}\sup_{0\leq s\leq b}|x^{\varphi,x_{1}}(s)-y^{\phi,y_{1}}(s)|\leq\epsilon, \quad\text{when}\quad\mathbf{E}\|\varphi-\phi\|_{\mathcal{B}}^{2}+\mathbf{E}|x_{1}-y_{1}|^{2}<\delta,
\end{equation}
where $y^{\phi,y_{1}}(t)$ is another solution of (\ref{eq1}) with initial value $(\phi,y_{1}).$
\end{defn}
\begin{thm}
Assume $21mM^{2}\sum_{t_{i}<t}p_{i}+21mM^{2}b\sum_{t_{i}<t}q_{i}<1$, the assumption of Theorem \ref{theorem} are satisfied and $g$ satisfied (H4), then the mild solution of (\ref{eq1}) is stable in mean square.
\end{thm}
\begin{proof}The proof is similar to the Theorem 18 in \cite{4}, we here give only the sketch of the proof.
 Using the same arguments as in Theorem \ref{theorem}, we get for all $0\leq t\leq b$ ,
 \begin{eqnarray*}
  &&\mathbf{E}\sup_{0\leq s\leq t}|x^{\varphi,x_{1}}(s)-y^{\phi,y_{1}}(s)|^{2}\\ &\leq& 3\mathbf{E}|S_{q}(t)\varphi(0)-S_{q}(t)\phi(0)|^{2}+3\mathbf{E}|\int_{0}^{t}S_{q}(s)(|x_{1}-y_{1}|_{H}+\|\varphi-\phi\|_{\mathcal{B}})\mathrm{d}s|^{2}\\
   && +3 (7mM^{2}\sum_{t_{i}<t}p_{i}+7mM^{2}b\sum_{t_{i}<t}q_{i})\mathbf{E}\sup_{0\leq s\leq t}|x^{\varphi,x_{1}}(s)-y^{\phi,y_{1}}(s)|^{2}\\
&&+3(14m^{2}M^{2}+7M^{2}b+7M_{b}^{2}\frac{b^{2q-1}}{2q-1}+7M_{b}^{2}b^{2q-2})\\&&\quad\times\int_{0}^{t}\tilde{\kappa}(\mathbf{E}\sup_{0\leq r\leq s}|x^{\varphi,x_{1}}(s)-y^{\phi,y_{1}}(s)|^{2})\mathrm{d}s.
 \end{eqnarray*}
 Then we get
 \begin{eqnarray*}
 &&\mathbf{E}\sup_{0\leq s\leq t}|x^{\varphi,x_{1}}(s)-y^{\phi,y_{1}}(s)|\\ &\leq&\frac{\nu}{\Lambda}(|x_{1}-y_{1}|^{2}+\|\varphi-\phi\|_{\mathcal{B}})+\frac{\tilde{\nu}}{\Lambda}\int_{0}^{t}\tilde{\kappa}(\mathbf{E}\sup_{0\leq r\leq s}|x^{\varphi,x_{1}}(s)-y^{\phi,y_{1}}(s)|^{2})\mathrm{d}s,
 \end{eqnarray*}
 where $\nu=\max\{6M^{2}b^{2},6M^{2}b^{2}K_{1}+3M^{2}L^{2}\}$ , $\tilde{\nu}=3(14m^{2}M^{2}+7M^{2}b+7M_{b}^{2}\frac{b^{2q-1}}{2q-1}+7M_{b}^{2}b^{2q-2}),$ and $\Lambda=1-3(7mM^{2}\sum_{t_{i}<t}p_{i}+7mM^{2}b\sum_{t_{i}<t}q_{i})$. Since the function $\tilde{\kappa}(u)$ is defined in (\ref{kk}) which has  the property as in Lemma \ref{B}. So for any $\epsilon>0$, letting $\epsilon_{1}=\frac{1}{2}\epsilon$, we have $\lim_{s\rightarrow0}\int_{s}^{\epsilon_{1}}\frac{1}{\tilde{\kappa}(u)}\mathrm{d}u=\infty.$ There exists a positive constant $\delta<\epsilon_{1}$ such that $\int_{\delta}^{\epsilon_{1}}\frac{1}{\tilde{\kappa}(u)}\mathrm{d}u\geq T.$ Let $u_{0}=\frac{\nu}{\Lambda}(|x_{1}-y_{1}|^{2}+\|\varphi-\phi\|_{\mathcal{B}}),$ $u(t)=\mathbf{E}\sup_{0\leq s\leq t}|x^{\varphi,x_{1}}(s)-y^{\phi,y_{1}}(s)|,$ $v(t)=1$. When $u_{0}\leq \delta\leq \epsilon_{1}$, Lemma \ref{cc} shows that $\int_{u_{0}}^{\epsilon_{1}}\frac{1}{\tilde{\kappa}(u)}\mathrm{d}u\geq \int_{\delta}^{\epsilon_{1}}\frac{1}{\tilde{\kappa}(u)}\mathrm{d}u\geq T=\int_{0}^{b}v(s)\mathrm{d}s.$ So for any $t\in[0,b],$ the estimate $u(t)\leq \epsilon_{1}\leq\epsilon$ holds. This completes the proof of the theorem.\end{proof}
\section{Application}
In this section, we present an example to show the obtained results.
Let $H=L^{2}([0,\pi])$, let's consider the following initial problem.
\begin{eqnarray}\label{eq2}
\left\{ \begin{array}{ll}
D_{t}^{\alpha}[u'_{t}(t,x)-\int_{-\infty}^{t}\int_{0}^{\pi}h(s-t,\eta,x)u(s,\eta)\mathrm{d}\eta\mathrm{d}s]\\=\frac{\partial^{2}}{\partial x^{2}}u(t,x)+f(t,\int_{-\infty}^{0}p_{0}(s)u(t-s,x)\mathrm{d}s)+\sigma(t,\int_{-\infty}^{0}q_{0}(s)u(t-s,x))\frac{\mathrm{d}w(t)}{\mathrm{d}t},\\
\qquad t\in J,\, t\neq t_{i};\\
\triangle u(t_{i},x)=\int_{-\infty}^{t_{i}}\bar{p}_{i}(s-t_{i})u(s,x)\mathrm{d}s,\,\\
\triangle u'(t_{i},x)=\int_{-\infty}^{t_{i}}\bar{q}_{i}(s-t_{i})\frac{u(s,x)}{1+|u(s,x)|}\mathrm{d}s,\,i=1,2,\ldots,m;\\
u(t,0)=u(t,\pi)=0,\,t\in[0,1];\quad u(\theta,x)=\varphi(\theta,x),\,\theta\in(-\infty,0],x\in[0,\pi];\\
\frac{\partial}{\partial t}u(0,x)=z(x),x\in[0,\pi].
\end{array} \right.
\end{eqnarray}
We present the abstract phase space $\mathcal{B}$ as,
  \begin{eqnarray*}
  &&\mathcal{B}=\{\psi:(-\infty,0]\rightarrow H:(E|\psi(\theta)|^{2})^{1/2} \,\textrm{is a bounded and measurable}\\ \\&&\textrm{function on }\,[-a,0]\,\textrm{and} \int_{-\infty}^{0}\rho(s)\sup_{s\leq\theta\leq0}(E|\psi(\theta)|^{2})^{1/2}\mathrm{d}s<\infty \}
  \end{eqnarray*}
  where $\rho:(-\infty, 0]\rightarrow(0,\infty)$ is a continuous function with $l=\int_{-\infty}^{0}\rho(t)\mathrm{d}t<\infty.$
  Let $\|\psi\|_{\mathcal{B}}=\int_{-\infty}^{0}\rho(s)\sup_{s\leq\theta\leq0}(E|\psi(\theta)|^{2})^{1/2}\mathrm{d}s,$ for all $\psi\in\mathcal{B}$, then $(\mathcal{B},\|\psi\|_{\mathcal{B}})$ is a Banach space (see \cite{6}).
 Obviously the operator $A:H\rightarrow H$ by $A=\frac{\partial^{2}}{\partial x^{2}}$ with domain $D(A)=\{z\in H: z(0)=z(\pi)=0\}$ is the infinitesimal generator of a strongly continuously cosine family\cite{15}. So $\|S_{\gamma}(t)\|_{\mathcal{L}(G,H)}\leq M$, $\|T_{\gamma}(t)\|_{\mathcal{L}(G,H)}\leq M_{1}.$
We suppose that
\begin{enumerate}
  \item [(a)] $h(s,\eta,x)$, $\frac{\partial h(s,\eta,x)}{\partial x}$ are measurable, $h(s,\eta,0)=h(s,\eta,\pi)=0$ ,
  \begin{eqnarray*}
  L_{0}=\max\{[\int_{0}^{\pi}\int_{-\infty}^{0}\int_{0}^{\pi}\frac{1}{\rho(s)}(\frac{\partial^{k} h(s,\eta,x)}{\partial x^{k}})^{2}\mathrm{d}\eta\mathrm{d}s\mathrm{d}x]^{\frac{1}{2}}:k=0,1\}<\infty,
  \end{eqnarray*}
 and
 \begin{eqnarray*}
  |f(t,\xi)-f(t,\zeta)|&\leq& L_{1}\|\xi-\zeta\|_{\mathcal{B}},\,|\sigma(t,\xi)-\sigma(t,\zeta)|\\
   &\leq& L_{1}\|\xi-\zeta\|_{\mathcal{B}}, \, t\in[0,1],\,\xi,\zeta\in \mathcal{B},
 \end{eqnarray*}
   where $L_{1},$ $L_{2},$  are two positive constants.
\item [(b)] $\bar{p}_{i}(\theta)\in C(\mathbf{R},\mathbf{R}^{+ })$ and $p_{i}=(\int_{-\infty}^{0}\frac{\bar{p}_{i}^{2}(\theta)}{\rho(\theta)}\mathrm{d}\theta)^{\frac{1}{2}}<\infty,$ $i=1,2,\ldots,m.$
  \item [(c)] $\bar{q}_{i}(\theta)\in C(\mathbf{R},\mathbf{R}^{+ })$ and $q_{i}=(\int_{-\infty}^{0}\frac{\bar{q}_{i}^{2}(\theta)}{\rho(\theta)}\mathrm{d}\theta)^{\frac{1}{2}}<\infty,$ $i=1,2,\ldots,m.$
\end{enumerate}
Let $L=\max\{L_{0},L_{1},L_{2}\}$, since we can take $\kappa(\vartheta)=L\cdot\vartheta$ which is a concave function as in (H1), so (H1)-(H3) hold.

By computation we have
\begin{eqnarray*}
 &&|\int_{-\infty}^{t}\int_{0}^{\pi}h(s-t,\eta,x)u(s,\eta)\mathrm{d}\eta\mathrm{d}s-\int_{-\infty}^{t}\int_{0}^{\pi}h(s-t,\eta,x)\tilde{u}(s,\eta)\mathrm{d}\eta\mathrm{d}s|\\
 &&\leq L_{0}\|u(t,\eta)-\tilde{u}(t,\eta)\|_{\mathcal{B}}
\end{eqnarray*}
 and it is obviously $\|u_{t}\|_{\mathcal{B}}\leq l\sup_{0\leq s\leq t}(E \|u(s,x)\|^{2})^{1/2}+\|\varphi(\theta,x)\|_{\mathcal{B}}.$
According to Theorem \ref{theorem} we get the conclusion that when
\begin{equation*}
 \max\{7mM^{2}l\sum_{i=1}^{m}p_{i}+14mM^{2}l\sum_{i=1}^{m}q_{i},7mM^{2}\sum_{i=1}^{m}p_{i}+7mM^{2}\sum_{i=1}^{m}q_{i}\}<1,
\end{equation*}
 problem (\ref{eq2}) has a unique mild solution on $(-\infty,1].$

\section*{}
SHUFEN ZHAO

Department of Mathematics, Harbin Institute of Technology, Harbin 150001, PR China\\
Department of  Mathematics , Zhaotong University, Zhaotong 657000,PR China\\
E-mail address: 12b312003@hit.edu.cn

MINGHUI SONG

Department of Mathematics, Harbin Institute of Technology, Harbin 150001, PR China\\
E-mail address: songmh@hit.edu.cn

\end{document}